# Prime-Power Expansions and Boundary Corrections for the Secular Coefficients of Riemann-Zero Windows

Francesco Cordella

## Abstract

We study the secular coefficients of characteristic polynomials attached to finite windows of zeros of the Riemann zeta function, mapped onto the unit circle. Starting from the Riemann–von Mangoldt formula and Newton's identities, we derive a symmetric prime-power model for their generating function. The model carries an explicit boundary factor $(1 - w)^{-\Delta S}$, where $\Delta S$ is the variation of $S(t)$ across the window. The boundary correction is determined by the integer zero count of the window rather than by a fitted parameter.

For finite prime-power tails we prove a mean-square identity with an explicit diagonal term and an explicit error bound, using the Montgomery–Vaughan generalisation of Hilbert's inequality.

Numerical experiments near $T = 10^6$, on computed but uncertified zeros, show that prime powers up to $Q = 32$, the zero count of the window and the known smooth term reproduce the secular coefficients. The median relative residual energy is $1.3 \times 10^{-4}$ for $H = 181$ and $1.0 \times 10^{-5}$ for $H = 724$; the means are $2.1 \times 10^{-3}$ and $2.7 \times 10^{-3}$. The error is strongly asymmetric across windows. An exact convolution identity separates a residual, whose measured size does not depend on the boundary exponent, from its amplification by the boundary factor; the measured amplification follows the exponent over eight orders of magnitude. The diagonal prediction matches the measured power-sum residuals within 16% in the tested regimes.

Exact identities and finite Dirichlet-polynomial estimates are kept separate from the numerical reconstruction. Extending the mean-square estimate to the true $S(t)$ is open and appears tied to pair correlations of zeros. This note does not propose a route to the Riemann hypothesis.

## 1. Introduction

The object of study is the monic polynomial whose roots are the zeros of zeta in a window (T, T+H], placed on the unit circle by a linear map. Its low coefficients turn out to be computable from a few prime powers, one integer per window (the zero count N) and a smooth term known in closed form.

Let $\gamma$ run over the ordinates of zeros with $T < \gamma \le T+H$, and let N be their number. We set

$$\omega_k = \frac{2\pi k}{H}, \qquad p_k = \sum_{T<\gamma\le T+H} e^{i\omega_k(\gamma - T)}, \qquad C(w) = \prod_{T<\gamma\le T+H} \left(1 - e^{i\omega_1(\gamma-T)} w\right) = \sum_{k=0}^{N} c_k\, w^k.$$

The numbers $c_k$ are, up to the sign $(-1)^k$, the elementary symmetric functions of the roots: the secular coefficients.

The note has three layers, and each statement below carries one of these labels.

- **Proved.** Exact identities for power sums and for the symmetric prime-power model (Sections 2–3), with the closed form of 3.2–3.3 as the central result, and a mean-square identity for finite prime-power tails (Section 4). Proofs are short and rest on two quoted results: a lemma of Ford and Zaharescu and an inequality of Montgomery and Vaughan.

- **Numerical.** Reconstruction of the true coefficients near $T = 10^6$, in modulus and phase (Section 5).
- **Open.** The mean-square estimate for the true remainder $S - S_Q$, its propagation through the Newton exponential, and the error term of the expansion over integers in 3.4 (Sections 3.4 and 6).

The main numerical result is a *conditioned* reconstruction: the boundary exponent is fixed by the zero count N of the same window, through $\Delta S = (\vartheta(T+H) - \vartheta(T))/\pi - N$. No continuous parameter is fitted, but it is not a prediction from primes alone.

Originality is not established. The correspondence between coefficients and integers is not new: Section 1.1 lists the related work found in a web search and says what was not found. The commensurate case of 2.3 belongs to the classical theory of para-orthogonal polynomials on the unit circle.

## 1.1 Related work

The dictionary between zeros, primes and coefficients is known, mostly at a heuristic level. We did not find the exact finite-window identity with its boundary factor, nor the numerical reconstruction. The search used the web, abstracts and excerpts, and is not exhaustive.

- **Windows of zeros on the circle.** Coram and Diaconis wrapped blocks of consecutive zeros around the unit circle and tested their traces against random unitary matrices; Diaconis's survey reports about 1 190 blocks, each wrapped as 42 points on the circle, from 50 000 consecutive zeros starting at the $10^{20}$-th. No arithmetic content is extracted.
- **Local polynomial model.** This is the closest precedent for Sections 2.3 and 3.4. Tao approximates zeta near height t by a polynomial of degree N whose coefficients are block sums of $n^{-1/2-it}$, with a self-reciprocal relation between coefficients j and N − j. In a second post the traces are block sums of $\Lambda(n)\, n^{-1/2-it}$. Both statements are heuristic, for N growing slowly with T, so each block holds many integers. Cut-off effects at the ends are noted and set aside.
- **Secular coefficients.** Haake and coauthors study the secular coefficients of random unitary matrices through Newton's formulae. Najnudel, Paquette and Simm relate them to a holomorphic multiplicative chaos. Soundararajan and Zaman study the coefficients of $\exp(\Sigma X(k) z^k / \sqrt{k})$ as a model for random multiplicative functions. These are random counterparts of Section 3.4. In quantum chaos the same step is the pseudo-orbit expansion of spectral determinants of Berry and Keating.
- **Moments.** Conrey and Gamburd relate moments of partial sums of zeta on the critical line to moments of truncated characteristic polynomials of random unitary matrices. Heap and Lindqvist compute moments of sums of random multiplicative functions and of truncated characteristic polynomials. That is a correspondence of moments, not of individual coefficients.
- **Primes and zeros in one product.** Gonek, Hughes and Keating write zeta as a partial Euler product times a partial Hadamard product over nearby zeros.
- **Power sums over zeros.** Landau's formula and its uniform versions by Gonek and by Ford and Zaharescu. Section 4.1 is a corollary.

- **Variance of S over intervals.** Under RH, Fujii's formula for the second moment of S(t+h) − S(t), built on Goldston's, contains F(α) for α ≥ 1 (as quoted by Chan). Lugar, Milinovich and Quesada-Herrera prove, under RH and a conjecture of Chan on gaps between zeros, a conjecture of Berry for the number variance in the non-universal regime, where GUE statistics do not describe the zeros. Both bear on ΔS and on Open problem 1.

- **Macroscopic statistics.** Rodgers proves, under RH and for sufficiently smooth test functions, that the pair-correlation measure of the zeros agrees closely with the one conjectured by Bogomolny and Keating, which has troughs at macroscopic scale.

- **Zeros from primes and an index.** LeClair computes high zeros from a truncated Euler product and the index n of the zero, without evaluating zeta.

- **Fisher–Hartwig singularities.** The factor $(1 - w)^{-\Delta S}$ is a Fisher–Hartwig singularity with a root-type (modulus) and a jump-type (phase) component, as in the theory of Toeplitz determinants (Deift, Its and Krasovsky). This link is our remark; we did not find it in the zeta literature.

Not found in this search:

- the identity of 3.2–3.3 with the factor $(1 - w)^{-\Delta S}$ fixed by the zero count;

- the regime where single integers are resolved, with sinc leakage and counter-rotating terms;

- the reconstruction of Section 5 and the dependence of its error on ΔS.

## 2. Exact identities

The secular coefficients are determined by the power sums through Newton's identities, and the power sums are windowed Fourier coefficients of S(t). Both statements are exact.

### 2.1 Newton

$$C(w) = \exp\left(-\sum_{k\geq 1} \frac{p_k}{k}\, w^k\right), \qquad c_0 = 1, \quad c_k = -\frac{1}{k}\sum_{j=1}^{k} p_j\, c_{k-j}.$$

If two families of power sums differ by $\delta_k$, the generating functions differ by the factor $\exp(-\Sigma\, \delta_k w^k / k)$. This is how errors propagate in Section 6.

### 2.2 Power sums and S(t) (proved)

Let $S(t) = \pi^{-1} \arg \zeta(1/2 + it)$ and $N(t) = \vartheta(t)/\pi + 1 + S(t)$. Assume T and T+H are not ordinates of zeros. Stieltjes integration by parts, with $e^{i\omega_k H} = 1$, gives for every k ≥ 1

$$p_k = A_k - \Delta S - i\omega_k \int_0^H e^{i\omega_k u}\, S(T+u)\, du, \qquad \Delta S = S(T) - S(T+H),$$

$$A_k = \frac{1}{\pi}\int_T^{T+H} e^{i\omega_k(t-T)}\, \vartheta'(t)\, dt = -\frac{iH^2}{4\pi^2 kT} + O\left(\frac{H^3}{kT^2}\right) \qquad (H \leq T).$$

No hypothesis on the location of the zeros is used. The number ΔS is computable from the window alone: $\Delta S = (\vartheta(T+H) - \vartheta(T))/\pi - N$. The smooth term $A_k$ depends only on T and H. Its asymptotic, for T ≥ 2, follows from one integration by parts and $\vartheta'(t) = \tfrac{1}{2}\log(t/2\pi) + O(t^{-2})$.

### 2.3 The commensurate model (proved, classical)

Let $A(z) = \Sigma_{k \le m} a_k z^k$ with $a_0 \neq 0$ and $m < n/2$. Write $\bar{A}(w) = \Sigma_{k \le m} \mathrm{conj}(a_k)\, w^k$ for the polynomial with conjugated coefficients, and let $P(z) = A(z) + c\, z^n \bar{A}(1/z)$ with $|c| = 1$.

- If A has no zeros in the closed unit disk, the n zeros $z_1, \ldots, z_n$ of P are simple and lie on the unit circle.
- Their elementary symmetric functions satisfy $e_k = (-1)^k \mathrm{conj}(a_k / a_0)$ for $k \le m$. Equivalently, $\Pi(1 - z_j w) = \bar{A}(w) / \mathrm{conj}(a_0) + O(w^{n-m})$.
- Their power sums satisfy $p_k = -k\, [w^k] \log \bar{A}(w)$ for $1 \le k < n - m$.

A numerical check with five frequencies reproduced the amplitudes exactly and placed the zeros on the circle to $10^{-14}$.

### 2.4 Non-integer frequency: linear response (proved)

For $0 < \kappa < n/2$ and complex $\varepsilon$, let $x_0, \ldots, x_{n-1}$ be the zeros in $[0, 2\pi]$ of

$$F_\varepsilon(x) = 2\,\mathrm{Re}\big[e^{-inx/2}\left(1 + \varepsilon e^{i\kappa x}\right)\big].$$

For $\varepsilon = 0$ these are the n simple zeros $(2j+1)\pi/n$. By the implicit function theorem there are still exactly n zeros, simple and smooth in $\varepsilon$, for $|\varepsilon| < \varepsilon_0(n, \kappa)$. Let $e_k(\varepsilon)$ be the elementary symmetric functions of the points $e^{ix_j}$. Then for $1 \le k < n$ and $|\varepsilon| < \varepsilon_0$, with an O-constant depending on n, k and $\kappa$,

$$e_k(\varepsilon) = (-1)^{k-1}[\varepsilon D_n(k + \kappa) - \bar{\varepsilon} D_n(k - \kappa)] + O(|\varepsilon|^2), \qquad D_n(u) = e^{i\pi u}\, \frac{\sin \pi u}{n \sin(\pi u/n)}.$$

The proof perturbs the roots of $z^n + 1$ to first order and linearises Newton's identities. Since $|D_n(u)|^2$ equals $\mathrm{sinc}^2(\pi u)$ up to a factor $1 + O(u^2/n^2)$, a frequency that is not an integer leaks into neighbouring coefficients with a $\mathrm{sinc}^2$ profile. For $n = 346$ and $\kappa = 24.9$ the measured error divided by $|\varepsilon|^2$ was about 0.006.

## 3. Symmetric prime-power model

Replacing S by its truncated prime-power expansion inside the integral of 2.2 gives a model whose Newton exponential has a closed form: an analytic part depending on the primes, times the boundary factor $(1 - w)^{-\Delta S}$. This closed form, with $\Delta S$ fixed by the zero count, is the central result of the note (3.2–3.3).

### 3.1 Definitions

Let q run over prime powers and let

$$g_Q(\theta) = \sum_{q \le Q} b_q\, e^{-i\kappa_q \theta}, \qquad b_q = \frac{\Lambda(q)}{\log q}\, q^{-1/2 - iT}, \qquad \kappa_q = \frac{H \log q}{2\pi},$$

so that $g_Q$ is the truncated Dirichlet series of $\log \zeta$ on the window, and $S_Q = \mathrm{Im}\, g_Q / \pi$. We also set

$$D(u) = e^{i\pi u}\, \frac{\sin \pi u}{\pi u}, \qquad \hat{f}_k = \frac{1}{2\pi}\int_0^{2\pi} f(\theta)\, e^{ik\theta}\, d\theta, \qquad J = g_Q(0) - g_Q(2\pi), \qquad \Delta S_Q = \frac{\mathrm{Im} J}{\pi}.$$

The symmetric model for the power sums keeps rotating and counter-rotating terms:

$$M_k^{\text{sym}} = -\sum_{q\le Q} \kappa_q \left[b_q\, D(k-\kappa_q) + \overline{b_q}\, D(k+\kappa_q)\right] = -\Delta S_Q - i\omega_k \int_0^H e^{i\omega_k u}\, S_Q(T+u)\, du.$$

### 3.2 Closed form (proved; central identity)

Two elementary identities hold for k ≥ 1:

$$\frac{\kappa}{k} D(k-\kappa) = D(k-\kappa) + \frac{1-e^{-2\pi i\kappa}}{2\pi i\, k}, \qquad \frac{\kappa}{k} D(k+\kappa) = -D(k+\kappa) - \frac{1-e^{2\pi i\kappa}}{2\pi i\, k}.$$

Summing over q and exponentiating with Newton gives

$$-\frac{M_k^{\text{sym}}}{k} = \left(\widehat{g_Q - \overline{g_Q}}\right)_k + \frac{\Delta S_Q}{k}, \qquad \sum_{k\ge 0} \tilde{c}_k\, w^k = \exp\left[\sum_{k\ge 1} \left(\widehat{g_Q - \overline{g_Q}}\right)_k w^k\right] (1-w)^{-\Delta S_Q}.$$

The boundary factor depends only on the imaginary part of J. The dependence on the frequencies stays in the analytic factor.

### 3.3 The full conditioned predictor (definition)

Two ingredients of 2.2 are known exactly and are restored: the true ΔS, fixed by the integer N through $\Delta S = (\vartheta(T+H) - \vartheta(T))/\pi - N$, and the smooth term $A_k$. No parameter is fitted. The full predictor is

$$p_k^{(Q,\text{full})} = A_k + M_k^{\text{sym}} - (\Delta S - \Delta S_Q), \qquad \tilde{C}_Q^{\text{full}}(w)$$
$$= \exp\left[-\sum_{k\ge 1} \frac{A_k}{k} w^k\right] \exp\left[\sum_{k\ge 1} \left(\widehat{g_Q - \overline{g_Q}}\right)_k w^k\right] (1-w)^{-\Delta S}.$$

By 2.2 the error of the predicted power sums is exactly

$$p_k - p_k^{(Q,\text{full})} = -i\omega_k \int_0^H e^{i\omega_k u}\, U_Q(T+u)\, du, \qquad U_Q = S - S_Q.$$

All of the error sits in the approximation of S(t) by $S_Q$; the smooth part is treated exactly. The predictor without the smooth term, $p_k^{Q,\text{corr}} = M_k^{\text{sym}} - (\Delta S - \Delta S_Q)$, is kept in Section 5 as an ablation test.

### 3.4 Where the integers come from (model identities; error estimate open)

Dropping the counter-rotating terms gives a simpler model. Its structure is exact, but its approximation by a sum over integers is not yet proved. Three statements are exact.

**(a) Generating function.** With rotating terms only, the first identity of 3.2 gives

$$\sum_{k\ge 0} \tilde{c}_k^{\text{rot}}\, w^k = \exp\left[\sum_{k\ge 1} \left(\widehat{g_Q}\right)_k w^k\right] (1-w)^{-J/2\pi i}.$$

**(b) Euler product.** The exponential of $g_Q$ is an absolutely convergent Dirichlet series over Q-smooth integers,

$$e^{g_Q(\theta)} = \sum_{m\ Q\text{-smooth}} a_Q(m)\, m^{-1/2-iT}\, e^{-i\kappa_m \theta}, \qquad \kappa_m = \frac{H \log m}{2\pi}.$$

Here $a_Q$ is multiplicative, and for each prime p its values on powers of p are the Taylor coefficients of $\exp(\Sigma\, x^i / i)$, the sum running over i with $p^i \le Q$. So $a_Q(m) = 1$ when every prime-power component of m is at most Q, and $0 < a_Q(m) < 1$ otherwise; for example $a_{13}(16) = 3/4$. Composite integers appear because $\kappa_{mn} = \kappa_m + \kappa_n$.

**(c) Projection.** Let $P_+$ and $P_-$ be the projections onto the Fourier modes $k \ge 1$ and $k \le -1$, so that $g_Q = \hat{g}_0 + P_+ g_Q + P_- g_Q$, and put $w = e^{-i\theta}$. Then

$$\exp(P_+ g_Q) = e^{-\widehat{(g_Q)}_0}\, e^{g_Q}\, e^{-P_- g_Q},$$

and the coefficients in (a) are exactly the Fourier coefficients of the bounded function

$$e^{-\widehat{(g_Q)}_0} \left[ \sum_{m\ Q\text{-smooth}} a_Q(m)\, m^{-1/2-iT}\, e^{-i\kappa_m \theta} \right] \left(1 - e^{-i\theta}\right)^{-J/2\pi i}\, e^{-P_- g_Q(\theta)}.$$

**Not proved.** Dropping the last factor gives, for $\Delta S_Q < 2$, the expansion over integers

$$\tilde{c}_k^{\,\mathrm{rot}} = e^{-\widehat{(g_Q)}_0} \sum_{m\ Q\text{-smooth}} a_Q(m)\, m^{-1/2-iT}\, D_J(k - \kappa_m) + R_k,$$

where $D_J$ is the Fourier kernel of $(1 - e^{-i\theta})^{-J/2\pi i}$. Each isolated term has amplitude $a_Q(m)\, m^{-1/2}$, so intensity at most 1/m; observed peaks are sums of such terms and include interference.

The remainder $R_k$ is the effect of the factor $e^{-P_- g_Q} - 1$. Two bounds are available: $|\hat{g}_0| \le \Sigma|b_q| / (\pi\, \kappa_2)$, and the $L^2$ norm of $P_- g_Q$ is at most $\Sigma|b_q| / (\pi \sqrt{\kappa_2})$. Both are small for long windows. But $P_- g_Q$ has a logarithmic singularity at the junction $\theta = 0$, of strength $|J|/2\pi$, so it is not uniformly small. An $L^2$ bound on it does not bound individual coefficients after exponentiation. A coefficient-wise estimate of $R_k$ is open.

Numerical check, for Q = 13 and 32, H from 181 to 11 585 and three heights each. With the factor $e^{-P_- g_Q}$ included, identity (c) holds to $10^{-6}$. Without it, the largest $|R_k|$ lies between $1.3 \times 10^{-4}$ and $5.5 \times 10^{-3}$, against coefficients of size up to 0.7. In the same runs the supremum of $|P_- g_Q|$ lies between 0.8 and 4.0, and its $L^2$ norm between 0.003 and 0.027. This small sample shows no clean decay of $R_k$ in H. Without the junction factor the discrepancy was about 0.16 |J|.

## 3.5 Effective boundary exponent of the predictor (proved)

For k large compared with $\kappa_Q$ the singular behaviour of the predictor is governed by $\delta_Q = \Delta S - \Delta S_Q$, not by $\Delta S$. The model carries its own boundary factor, which partly cancels the explicit one.

Integrating by parts in the last formula of 3.1, with $S_Q$ smooth, gives for fixed Q and H

$$M_k^{\mathrm{sym}} = \int_0^H e^{i\omega_k u}\, S_Q'(T+u)\, du = O_Q\!\left(\frac{1}{k}\right), \qquad \left(\widehat{g_Q - \overline{g_Q}}\right)_k = -\frac{\Delta S_Q}{k} + O_Q\!\left(\frac{1}{k^2}\right).$$

So the analytic factor of 3.2 already contains $(1 - w)^{\Delta S_Q}$. With $\delta_Q = \Delta S - \Delta S_Q$, the predictor of 3.3 is exactly

$$\tilde{C}_Q^{\,\mathrm{full}}(w) = F_Q(w)\,(1 - w)^{-\delta_Q}, \qquad F_Q(w) = \exp\left[ -\sum_{k\ge1} \frac{A_k + M_k^{\mathrm{sym}}}{k}\, w^k \right].$$

For $\delta_Q > -1$ and $\delta_Q \ne 0$ the coefficients of the predictor satisfy $\tilde{c}_k \sim F_Q(1)\, k^{\delta_Q - 1} / \Gamma(\delta_Q)$ as $k \to \infty$. In particular the predictor belongs to the Hardy space $H^2$ exactly when $\delta_Q < ½$, and for $\delta_Q > ½$ the sum of $|\tilde{c}_k|^2$ up to K grows like $K^{2\delta_Q - 1}$.

Proof. Put $\varphi_k = -(A_k + M_k^{sym})/k$, so $|\varphi_k| \le C/k^2$ by 2.2 and the bound above. Sequences a with sup $(1+k)^2 |a_k|$ finite form an algebra under convolution: splitting the sum at $j = k/2$ shows that the convolution of two such sequences is again $O(k^{-2})$. Hence the Taylor coefficients $f_k$ of $F_Q = \exp(\Sigma \varphi_k w^k)$ are $O(k^{-2})$. So $F_Q$ is continuous on the closed unit disk, and $F_Q(1) \ne 0$ because $F_Q$ is an exponential.

Let $d_r = (\delta_Q)_r / r!$ be the coefficients of $(1 - w)^{-\delta Q}$, so $d_r \sim r^{\delta_Q - 1} / \Gamma(\delta_Q)$, and write $\tilde{c}_k = \Sigma_j f_j d_{k-j}$. For $j \le k/2$ one has $d_{k-j} = d_k (1 + O(j/k))$, so these terms give $d_k [F_Q(1) + O(\log k / k)]$. The terms with $j > k/2$ are at most max $|f_j|$ times $\Sigma_{r<k/2} |d_r|$, which is $O(k^{\delta_Q - 2})$ for $\delta_Q > 0$ and $O(k^{-2})$ for $\delta_Q < 0$. Both are $o(d_k)$ when $\delta_Q > -1$.

Numerical check, for Q = 32, H = 181, six windows and k up to 6 000. The exponent fitted on $2\,000 \le k \le 6\,000$ agrees with $\delta_Q - 1$ within 0.012. In the five windows with $\delta_Q > 0$, the amplitude agrees with $|F_Q(1)| / \Gamma(\delta_Q)$ within 0.8%.

Scope. The asymptotic regime is k much larger than $\kappa_Q$. The range reconstructed in Section 5 is $k \le K$, with K below $\kappa_Q$ for $Q \ge 32$. So this statement does not by itself describe the errors measured there; Section 5.7 tests it.

### 3.6 Amplification over a finite range (proved)

The error of the coefficients is an exact convolution of a residual sequence with the coefficients of $(1 - w)^{-\Delta S}$. This holds for every k, with no asymptotic assumption.

Let $r_k = p_k - p_k^{Q,full}$ be the residual of 3.3, and write the predictor as $G_Q(w) (1 - w)^{-\Delta S}$, where $G_Q$ is the product of the two exponentials in 3.3. By 2.1 the true generating function is the predictor times $\exp(-\Sigma r_k w^k / k)$. Define

$$B_Q(w) = G_Q(w) \left[ \exp\left( -\sum_{k\ge1} \frac{r_k}{k} w^k \right) - 1 \right] = \sum_{j\ge1} b_j w^j.$$

Then the difference of the two generating functions is $B_Q(w) (1 - w)^{-\Delta S}$, that is,

$$c_k - \tilde{c}_k = \sum_{j=1}^{k} \frac{(\Delta S)_{k-j}}{(k-j)!} b_j, \qquad (a)_r = a(a+1)\cdots(a+r-1), \quad (a)_0 = 1.$$

Young's inequality for convolutions gives, for every $K \ge 1$,

$$\left( \sum_{k=1}^{K} |c_k - \tilde{c}_k|^2 \right)^{1/2} \le L_K(\Delta S) \left( \sum_{j=1}^{K} |b_j|^2 \right)^{1/2}, \qquad L_K(a) = \sum_{r=0}^{K-1} \frac{|(a)_r|}{r!}.$$

For $a > 0$ all terms are positive and $L_K(a) = \Gamma(K + a) / (\Gamma(K)\, \Gamma(1 + a))$, which is $K^a / \Gamma(1 + a)$ up to a factor $1 + O(1/K)$. For $-1 \le a \le 0$ one has $L_K(a) \le 2$ for every K. So for negative exponents in this range the amplification is bounded independently of K. It need not be below 1: $a = -½$, K = 2 and $b = (1, -1)$ give an energy ratio 13/8.

Young's bound is not sharp. The sharp constant is the norm of the K × K lower-triangular Toeplitz matrix of the coefficients $(a)_r / r!$. For a = 1 this norm is $1 / (2 \sin(\pi / (4K + 2)))$, about $2K/\pi$, so Young's bound overestimates the energy by a factor near 2.5.

The lemma separates two effects: the size of the residual before the boundary factor, $B_K = \Sigma_{j\le K} |b_j|^2$, and its amplification. It does not bound $B_K$, which depends on $S - S_Q$ (Open problem 1).

Summation by parts gives a second exact form. Let $h_k(a) = (a+1)_{k-1} / (k-1)!$ be the partial sums of the coefficients of $(1 - w)^{-a}$, and $\Delta b_j = b_{j+1} - b_j$. Then, for every k,

$$c_k - \tilde{c}_k = b_1\, h_k(\Delta S) + \sum_{j=1}^{k-1} \Delta\, b_j\, h_{k-j}(\Delta S).$$

If $b_j$ is constant the second term vanishes, and the ratio of the error energy to $B_K$ is exactly

$$\Phi_K(a) = \frac{1}{K}\sum_{k=1}^{K} h_k\,(a)^2 \ \sim\ \frac{K^{2a}}{(2a+1)\,\Gamma(1+a)^2} \qquad \left(a > -\frac{1}{2},\ K \to \infty\right).$$

So for a constant residual the error energy grows like a power of K when a > 0 and decays when −½ < a < 0; at a = −½ the leading term is log K / (πK). For a residual that varies, the second term of the identity has to be estimated.

A reason to expect slow variation: $-r_k/k$ is 2πi times the k-th Fourier coefficient of $S - S_Q$ on the window. For k ≤ K below $\kappa_Q$ every omitted frequency lies above the range, so these coefficients are leakage terms. No estimate of $\Delta b_j$ is proved; Section 5.7 measures it.

## 4. Mean-square estimates

Two rigorous statements are available: a pointwise formula for the power sums with a large error term, and a mean-square formula for the contribution of a finite tail of prime powers.

### 4.1 Windowed Landau formula (proved, from a lemma of Ford and Zaharescu)

Let k ≥ 1, $x = e^{\omega_k}$, q the prime power nearest to x, and $\kappa_q = H \log q / 2\pi$. Let T > 1 and assume all zeros with T < γ ≤ T+H lie on the critical line. Then, uniformly in k, T and H,

$$p_k = -\frac{H\,\Lambda(q)}{2\pi\sqrt{x}}\ q^{-iT}\, e^{i\pi(k-\kappa_q)}\, \frac{\sin\pi\left(k-\kappa_q\right)}{\pi\left(k-\kappa_q\right)} + E_k, \qquad E_k \ll \sqrt{x}\log^2\left(2x(T+H)\right) + \frac{\log\left(2(T+H)\right)}{\sqrt{x}\,\log x}.$$

The proof applies Lemma 1 of Ford–Zaharescu at T+H and at T, subtracts, and uses $x^\rho = \sqrt{x} \cdot x^{i\gamma}$ on the line. The hypothesis holds unconditionally for T+H ≤ 3 000 175 332 800 by Platt–Trudgian.

The bound is uniform but useful only in part of the range. For fixed q and k within a bounded distance of $\kappa_q$, the main term has size H Λ(q) / (2π√q) and the error bound has size √q $\log^2$(qT). So the formula is an asymptotic when H / (q $\log^2$ T) → ∞. For k small compared with $\kappa_2$ the second error term, of size H log T / k, exceeds the main term.

Where it applies, the statement gives the sinc profile and the peak height $\Lambda(q)^2 / (q\,\tau^2)$ for normalised power sums, with τ = log(T/2π). For short windows the bound is far from sharp: at T = $10^6$ and N = 346 the main term at q = 2 is 14, while the shape of the bound is about 470.

### 4.2 Finite prime-power tails (proved)

For Q < R let $U_{Q,R}(t) = -\pi^{-1} \sum_{Q<q\le R} \Lambda(q) \sin(t \log q) / (\sqrt{q} \log q)$, and let $R_{k,Q,R}(T)$ be its contribution to the error of 3.3. With $I_H(v) = (e^{ivH} - 1)/(iv)$,

$$R_{k,Q,R}(T) = \frac{\omega_k}{2\pi} \sum_{Q<q\le R} \frac{\Lambda(q)}{\sqrt{q}\,\log q} \left[e^{iT\log q}\, I_H(\omega_k + \log q) - e^{-iT\log q}\, I_H(\omega_k - \log q)\right].$$

This is one exponential sum in T with distinct real frequencies ± log q, two for each prime power q in (Q, R]. Rotating and counter-rotating terms are combined into a single exponential sum, so all cross terms are bounded together by the inequality below. The Montgomery–Vaughan generalisation of Hilbert's inequality states that, for distinct real $\lambda_r$ with $\delta_r$ the distance from $\lambda_r$ to the nearest other frequency,

$$\left| \int_X^{2X} \left| \sum_r a_r \, e^{i\lambda_r T} \right|^2 dT - X \sum_r |a_r|^2 \right| \le 3\pi \sum_r \frac{|a_r|^2}{\delta_r}.$$

The original statement is for the interval [0, X]; replacing $a_r$ by $a_r \, e^{i\lambda_r X}$ moves it to [X, 2X]. The frequencies here are well separated. Two prime powers q < q' give log q' − log q ≥ log(1 + 1/q) ≥ 1/(q+1), and a positive and a negative frequency differ by at least 2 log 2. So δ ≥ 1/(q+1) for both frequencies ± log q, and 1/δ ≤ R + 1.

**Theorem.** For every k ≥ 1, every 1 ≤ Q < R and every X > 0,

$$\left| \frac{1}{X} \int_X^{2X} \left| R_{k,Q,R}(T) \right|^2 dT - \mathcal{V}_{k,Q,R} \right| \le \frac{3\pi(R+1)}{X} \mathcal{V}_{k,Q,R},$$

$$\mathcal{V}_{k,Q,R} = \frac{\omega_k^2}{4\pi^2} \sum_{Q<q\le R} \frac{\Lambda(q)^2}{q \log^2 q} \left( |I_H(\omega_k + \log q)|^2 + |I_H(\omega_k - \log q)|^2 \right).$$

No condition on $\omega_k$ is needed. The statement is for a finite polynomial and says something only when R is small compared with X. It must not be applied to R = ∞.

A direct evaluation of the integral in nine cases (H = 181 and 724, R up to 3 000, X from $10^4$ to $10^6$) gave relative deviations from V between $10^{-5}$ and $6 \times 10^{-3}$, always below the bound.

Assume now $\omega_k \le \frac{1}{2} \log Q$, so that $|\omega_k - \log q| \ge \frac{1}{2} \log q$ for q > Q. Since $|I_H(v)| \le 2/|v|$ and $\Sigma_{p>Q} \, 1/(p \log^2 p) \ll 1/\log^2 Q$,

$$\frac{1}{X} \int_X^{2X} \left| R_{k,Q,R}(T) \right|^2 dT \ll \left(1 + \frac{R}{X}\right) \frac{\omega_k^2}{\log^2 Q}.$$

In the more restrictive regime $\omega_k / \log Q \to 0$, and if $\sin^2$ averages to ½ over the tail, the root mean square is close to $\omega_k / (\pi\sqrt{2} \log Q) \approx 0.225 \, \omega_k / \log Q$. This constant is a heuristic guide. It is not uniform when $\omega_k$ approaches log Q.

# 5. Numerical reconstruction

Near T = $10^6$ the full predictor of 3.3 with Q = 32 reproduces the true secular coefficients, in modulus and phase, with median relative residual energy $1.3 \times 10^{-4}$ for H = 181 and $1.0 \times 10^{-5}$ for H = 724. All results in this section are numerical. Values in square brackets are 95% moving-block bootstrap intervals over windows.

## 5.1 Data

We use 100 215 computed zeros with $10^6 < \gamma < 1.0525 \times 10^6$. Windows are fixed and adjacent: $T_j = \gamma_1 - 0.01 + jH$, where $\gamma_1$ = 1 000 000.5841 is the first zero above $10^6$, with j = 0, …, 288 for H = 181 and j = 0, …, 71 for H = 724. Residual energy means $\Sigma|c_k - \text{predicted}|^2 / \Sigma|c_k|^2$ over 1 ≤ k ≤ K.

| H | zeros per window | windows | K |
|---|---|---|---|
| 181 | 346 | 289 | 90 |
| 724 | 1 381 | 72 | 360 |

## 5.2 The pointwise formula of 4.1

The true error is of order 1 and does not grow with H, three orders of magnitude below the proved bound. The model here uses all prime powers up to 32.

| H | zeros per window | main term at q = 2 | rms of $E_k$ | max of $E_k$ |
|---|---|---|---|---|
| 181 | 346 | 14 | 0.60 | 2.0 |
| 724 | 1 381 | 55 | 0.73 | 1.6 |
| 2 896 | 5 522 | 150 | 0.52 | 1.1 |
| 11 585 | 22 096 | 902 | 0.38 | 0.8 |

Around q = 2 with H = 724 the mean of $|p_k|^2$ follows the $\mathrm{sinc}^2$ profile over three orders of magnitude: 3018 at the peak (predicted 3014), 67.8 one step away (68.7), 2.7 five steps away (3.1).

## 5.3 Where the discrepancy sits

The gap between true power sums and the symmetric model is almost entirely a constant in k: the difference between the true jump of S and its prime-power approximation, which converges slowly in Q (H = 181).

| Q | rms of $p_k - M_k^{sym}$ | rms of $\Delta S - \Delta S_Q$ |
|---|---|---|
| 32 | 0.52 | 0.46 |
| 3 000 | 0.36 | 0.36 |
| 300 000 | 0.31 | 0.31 |

## 5.4 Comparison of predictors

Median residual energy for three versions of the symmetric model. The first takes ΔS from primes. The second takes ΔS from the zero count but omits $A_k$ (ablation). The third is the full predictor.

| H, Q | ΔS from primes | zero count, no $A_k$ | full predictor |
|---|---|---|---|
| 181, 32 | 0.17 | $1.4 \times 10^{-4}$ | $1.3 \times 10^{-4}$ |
| 181, 3 000 | 0.13 | $1.4 \times 10^{-5}$ | $1.3 \times 10^{-5}$ |
| 181, 300 000 | 0.07 | $3.6 \times 10^{-6}$ | $2.6 \times 10^{-6}$ |
| 724, 32 | 0.23 | $2.4 \times 10^{-4}$ | $1.0 \times 10^{-5}$ |
| 724, 3 000 | 0.16 | $2.4 \times 10^{-4}$ | $1.6 \times 10^{-6}$ |
| 724, 300 000 | 0.08 | $2.3 \times 10^{-4}$ | $9.9 \times 10^{-7}$ |

The zero count is essential. Altering ΔS by 0.05 in one window (H = 724, ΔS = −0.02) raised the residual of the full predictor from $1.1 \times 10^{-5}$ to $5 \times 10^{-3}$; by 0.2, to 0.11. The smooth term matters for the longer window: without it the median for H = 724 stays near $2.4 \times 10^{-4}$ for every Q. With rotating terms only and a complex junction factor, the median for H = 181 and Q = 32 is 0.16.

### 5.5 Full predictor: distribution over windows

| H, Q | median [95% interval] | mean [95% interval] | 95th percentile | maximum |
|---|---|---|---|---|
| 181, 32 | $1.3 \times 10^{-4}$ [0.7, 2.5] | $2.1 \times 10^{-3}$ [1.5, 2.8] | $1.1 \times 10^{-2}$ | $5.3 \times 10^{-2}$ |
| 181, 3 000 | $1.3 \times 10^{-5}$ [0.8, 2.6] | $4.4 \times 10^{-4}$ [3.3, 5.8] | $2.9 \times 10^{-3}$ | $7.7 \times 10^{-3}$ |
| 724, 32 | $1.0 \times 10^{-5}$ [0.5, 8.0] | $2.7 \times 10^{-3}$ [1.2, 4.9] | $1.2 \times 10^{-2}$ | $4.8 \times 10^{-2}$ |
| 724, 3 000 | $1.6 \times 10^{-6}$ [0.6, 27] | $6.0 \times 10^{-4}$ [2.6, 11] | $3.1 \times 10^{-3}$ | $1.1 \times 10^{-2}$ |

Interval endpoints are in units of the power of ten of the same cell.

The medians for H = 724 rest on 72 windows, and their intervals span more than a decade. At the positions of the integers m = 2, 3, 4, 6, 10, 12 (Q = 300 000), the median ratio of predicted to true modulus is 1.000. The rms phase difference is 0.003–0.018 rad for H = 181 and 0.006–0.015 rad for H = 724. In the window of Figure 1 (Q = 32), the coefficient at m = 6 is 0.1032 − 0.2749i against a predicted 0.1007 − 0.2751i.

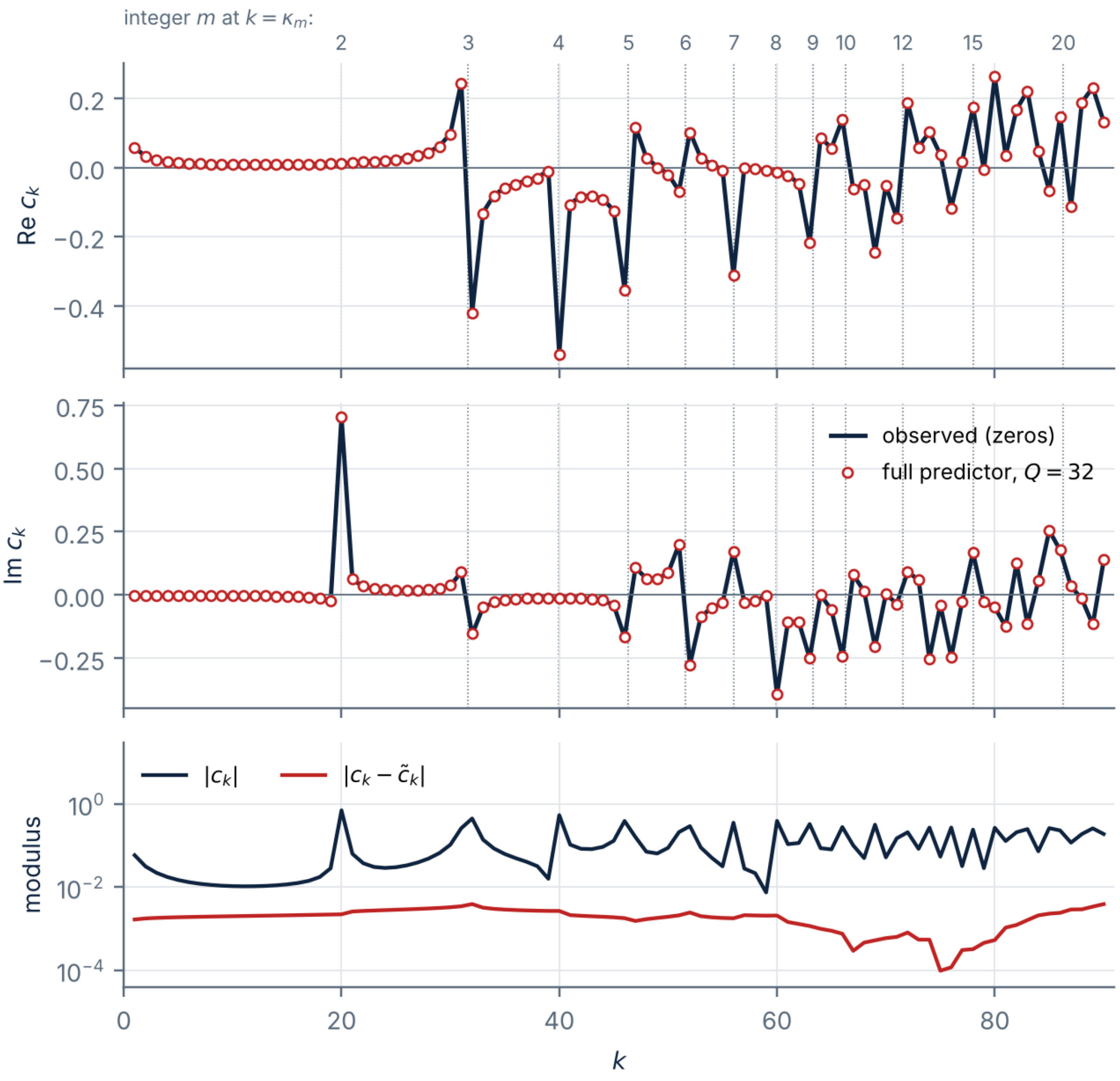

*Figure 1.* Observed and predicted secular coefficients $c_k$ for one window of 345 zeros (H = 181, T = 1 000 543.57, Q = 32, ΔS = 0.058, residual energy $1.4 \times 10^{-4}$). Dotted lines mark $k = \kappa_m$ for the integers m shown on top. Bottom panel: $|c_k|$ and the error of the full predictor.

### 5.6 Residual power sums against the diagonal term of 4.2

The residual $p_k - p_k^{Q,full}$ grows linearly with k. Its rms was compared with $\sqrt{V}$ computed with the sum over prime powers extended formally to infinity. The table gives the range of measured / predicted over all k from 1 to K. Figure 2 shows the same comparison.

| Window and range of k | Q = 32 | Q = 300 | Q = 3 000 |
|---|---|---|---|
| H = 181, k = 1 to 90 | 0.94–0.96 | 0.91–0.93 | 0.96–0.98 |
| H = 724, k = 1 to 360 | 0.91–0.95 | 0.96–0.99 | 1.14–1.16 |

For H = 181 and k = 90 the measured values are 0.425 (Q = 32) and 0.096 (Q = 3 000), against 0.442 and 0.098. Truncating the sum at R = T gives 0.084 for Q = 3 000, a worse fit. The ratio is nearly constant in k: the diagonal term gives the shape, and misses a constant factor between 0.91 and 1.16.

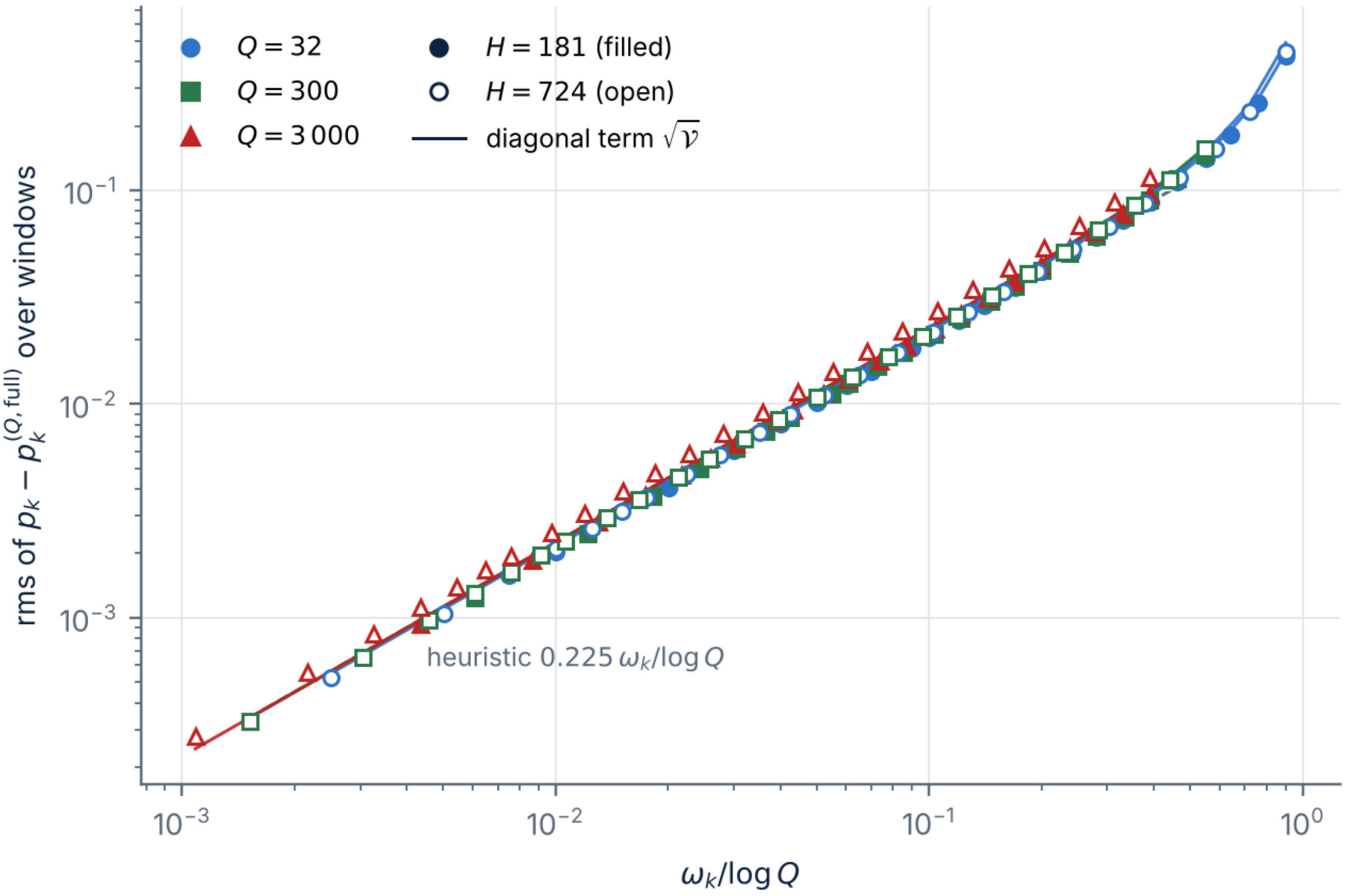

*Figure 2.* Root mean square over windows of the residual power sums $p_k - p_k^{Q,full}$, against $\omega_k$ / log Q. Solid lines: the diagonal term $\sqrt{V}$ of 4.2, with the prime sum extended formally to infinity. Dashed line: the heuristic 0.225 $\omega_k$ / log Q.

Ablation. Without $A_k$ the residual has a floor at small k equal to $|A_k| \approx H^2/(4\pi^2 kT)$. For H = 724, k = 4 and Q = 3 000 the rms is 0.0034 without $A_k$, against 0.0033 for the floor, and 0.0011 with $A_k$. For H = 724 this floor alone holds the median residual energy near $2.4 \times 10^{-4}$ for every Q (Section 5.4).

### 5.7 Anomalous windows

The mean residual energy is 15 to 400 times the median, and the excess is strongly associated with the sign and size of ΔS. The Spearman correlation between the logarithm of the residual energy and the signed ΔS is 0.87 [0.83, 0.89] for H = 181 and 0.94 [0.91, 0.96] for H = 724 (Q = 32). The association persists once the floor due to $A_k$ is removed, and for H = 724 it then spans seven orders of magnitude. Numerical conditioning is ruled out below. Figure 3 shows all windows.

| ΔS (H = 181, Q = 32) | windows | median residual energy | median of $\Sigma\|c_k\|^2$ |
|---|---|---|---|
| below −0.5 | 69 | $2.8 \times 10^{-6}$ | 4.6 |
| −0.5 to 0 | 62 | $2.6 \times 10^{-5}$ | 3.0 |
| 0 to 0.5 | 111 | $7.1 \times 10^{-4}$ | 3.1 |
| 0.5 to 0.8 | 22 | $5.8 \times 10^{-3}$ | 9.6 |
| above 0.8 | 25 | $4.7 \times 10^{-3}$ | 262 |

| ΔS (H = 724, Q = 32) | windows | median residual energy | median of $\Sigma\|c_k\|^2$ |
|---|---|---|---|
| below −0.5 | 19 | $1.2 \times 10^{-7}$ | 5.0 |
| −0.5 to 0 | 19 | $3.6 \times 10^{-6}$ | 2.9 |
| 0 to 0.5 | 16 | $3.1 \times 10^{-4}$ | 4.0 |
| 0.5 to 0.8 | 10 | $2.5 \times 10^{-3}$ | 32 |
| above 0.8 | 8 | $4.4 \times 10^{-3}$ | 1 127 |

Windows with ΔS > 0.5 are 16% of the sample and carry 61% [43, 77] of the total residual energy for H = 181; for H = 724 they are 25% and carry 94% [82, 99]. A mechanism consistent with this is the slow decay of the boundary factor. For ΔS not a non-positive integer,

$$[w^k]\,(1-w)^{-\Delta S} = \frac{\Gamma(k+\Delta S)}{\Gamma(\Delta S)\,\Gamma(k+1)} \sim \frac{k^{\Delta S-1}}{\Gamma(\Delta S)} \qquad (k \to \infty).$$

For ΔS near 1 these coefficients do not decay, and errors concentrated at the junction w = 1 are amplified; for ΔS < 0 they decay faster than 1/k. This is not a proof of the cause. By 3.5 the predictor itself behaves, for k much larger than $\kappa_Q$, like $(1-w)^{-\delta_Q}$ with $\delta_Q = \Delta S - \Delta S_Q$, so the exponent that matters for k ≤ K has to be settled by the data (below). The factor is a Fisher–Hartwig singularity with a root-type and a jump-type component (Section 1.1); its amplifying effect does not depend on the zeta function.

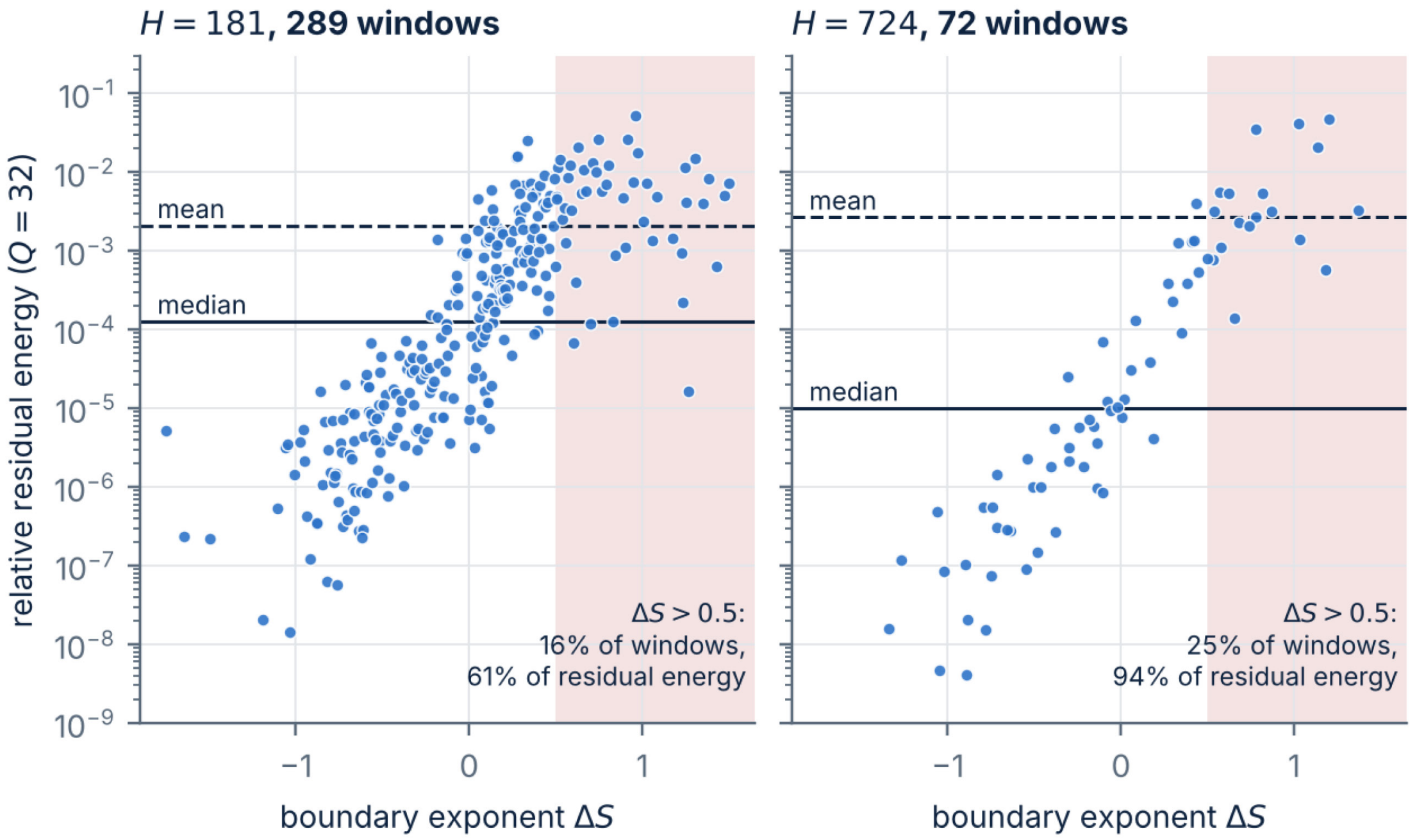

*Figure 3.* Relative residual energy of the full predictor (Q = 32) against the boundary exponent ΔS, one point per window. Shaded region: ΔS > 0.5.

Which exponent. Two candidates are the true jump $\Delta S$ and the effective exponent $\delta_Q = \Delta S - \Delta S_Q$ of 3.5. In the range $k \le K$ the data favour $\Delta S$. The table gives Spearman correlations of the logarithm of the residual energy, and partial rank correlations.

| H, Q | with ΔS | with $\delta_Q$ | $\delta_Q$ given ΔS | ΔS given $\delta_Q$ |
|---|---|---|---|---|
| 181, 32 | 0.87 | 0.67 | −0.03 | 0.75 |
| 181, 3 000 | 0.83 | 0.51 | −0.03 | 0.75 |
| 724, 32 | 0.94 | 0.80 | 0.24 | 0.85 |
| 724, 3 000 | 0.86 | 0.66 | −0.04 | 0.73 |

Once $\Delta S$ is known, $\delta_Q$ carries almost no further information; the converse does not hold. The threshold $\delta_Q = ½$ also separates the windows less sharply: for H = 181 and Q = 32, windows with $\delta_Q > ½$ carry 40% of the residual energy, against 61% for $\Delta S > ½$. A plausible reason is the scope remark of 3.5: K lies below $\kappa_Q$, where the cancellation between the two boundary factors has not yet set in. This has not been derived.

The association is not produced by the normalisation. For Q = 32 the absolute error $\Sigma|c_k - \tilde{c}_k|^2$ has Spearman correlation 0.90 (H = 181) and 0.95 (H = 724) with $\Delta S$, while the signal $\Sigma|c_k|^2$ has 0.21 and 0.40. Windows with $\Delta S > ½$ carry more than 99% of the total absolute error.

Residual and amplification. By 3.6 the absolute error is $B_K$ times an amplification ratio. The two behave very differently. The table gives Spearman correlations with $\Delta S$ for Q = 32.

| Quantity | H = 181 | H = 724 |
|---|---|---|
| $B_K$, residual before the boundary factor | 0.02 | −0.12 |
| $\Sigma|c_k - \tilde{c}_k|^2$, absolute error | 0.90 | 0.95 |
| their ratio, the amplification | 0.99 | 0.99 |

For Q = 3 000 the three values are −0.01, 0.85, 0.99 (H = 181) and −0.15, 0.87, 0.98 (H = 724). The amplification accounts for 85% of the variance of the logarithm of the absolute error for H = 181 and for 98% for H = 724.

The residual before amplification has the size expected from Section 4.2. The mean of $\Sigma|r_k/k|^2$ over windows is between 0.86 and 1.32 times $\Sigma\, V_k/k^2$, and the mean of $B_K$ is 1.4 to 1.6 times larger.

The amplification ratio spans eight to nine orders of magnitude: from $4 \times 10^{-4}$ to $5 \times 10^4$ for H = 181, and from $3 \times 10^{-4}$ to $4 \times 10^5$ for H = 724. The inequality of 3.6 holds in every window. The ratio to Young's bound is close to 1 only near $\Delta S = 0$, where the bound is itself close to 1. For $\Delta S > ½$ its median is 0.24 and its maximum 0.45 (H = 181), or 0.23 and 0.38 (H = 724). Relative to the sharp constant of 3.6 the medians are 0.54 and 0.53.

The law $\Phi_K(\Delta S)$, exact for a constant residual, reproduces the ratio. Its median quotient with the data is 0.89 and 0.88, 90% of the windows lie between 0.24 and 3.3, and the correlation of the logarithms is 0.99.

The residual sequences are smooth but not constant. The relative variation $F_K = \Sigma|b_j - \bar{b}|^2 / \Sigma|b_j|^2$, with $\bar{b}$ the mean of $b_j$, has median 0.28 for H = 181 and 0.31 for H = 724; the ratio $\Sigma|\Delta b_j|^2 / \Sigma|b_j|^2$ has median 0.012 and 0.002. For $\Delta S > ½$ the error differs from the flat term $\bar{b}\, h_k(\Delta S)$ by 7% of its energy (H = 181) and 14% (H = 724), in the median. For $\Delta S < -½$ it does not: there the variation term of 3.6 dominates, and $\Phi_K$ gives the order of magnitude only.

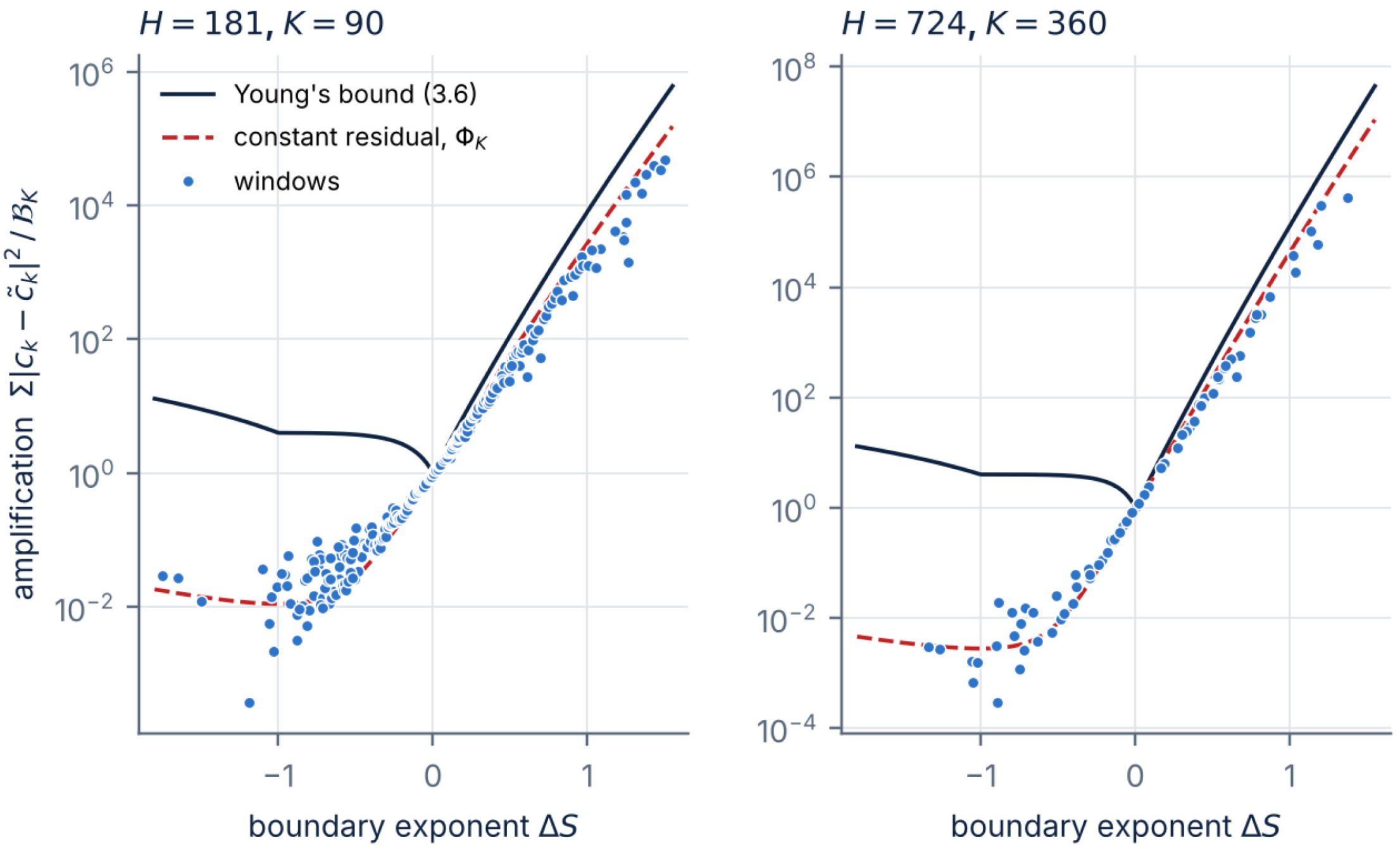

*Figure 4.* Amplification ratio $\Sigma|c_k - \tilde{c}_k|^2 / B_K$ against ΔS for Q = 32, one point per window. Solid line: Young's bound $L_K(\Delta S)^2$ of 3.6. Dashed line: $\Phi_K(\Delta S)$, exact for a constant residual.

So in these data the anomalous windows are explained by amplification. The residual is of ordinary size in every window, and the boundary factor multiplies it by a factor that depends almost only on ΔS and K. For negative ΔS the same convolution attenuates the residual sequences observed here, which is why those windows are reconstructed best. This attenuation is numerical; it is not a general property (Section 3.6). The rigorous part is the lemma of 3.6; the size of $B_K$ and the agreement with $\Phi_K$ are numerical.

Other candidates were weak. The minimum normalised gap correlates at 0.05–0.26 and |Z| at the endpoints at about 0. The distance from the endpoints to the nearest zero correlates at 0.46–0.65, which the sign of ΔS explains: S is largest just after a zero. Newton's recurrence in double precision agrees with a 60-digit computation to $5 \times 10^{-16}$ on the worst window. The zero count of every window was confirmed with an independent implementation (Appendix A.1), which excludes miscounted zeros as a source of the association.

## 6. Discussion and open problems

The proved part covers the model and finite prime-power tails; the accuracy on true zeros is numerical. The gap between the two is one well-defined estimate.

### 6.1 What the reconstruction is

The predictor uses a few prime powers, the smooth term $A_k$ and one integer per window, the zero count N, which fixes ΔS. It is a reconstruction conditioned on that count. With ΔS taken from primes the median residual energy stays between 0.23 and 0.07 for Q from 32 to 300 000, because prime sums approximate the jump of S slowly.

### 6.2 Relation to the mean square of S(t)

Assuming RH, Goldston's formula for the mean square of S(t) contains, beyond $(T/2\pi^2) \log\log T$, a term of order T that involves Montgomery's pair-correlation function $F(\alpha)$ for $\alpha \geq 1$ (statement as quoted here). This indicates that correlations between zeros enter the variance of S at order 1. It does not by itself give the variance of the part of S not represented by primes up to T: that also requires the covariance between S and the subtracted Dirichlet polynomial.

In our data the diagonal term with the prime sum extended formally beyond T fits better than the sum truncated at T (Section 5.6). This is consistent with a zero-dependent contribution of order 1, and it is numerical evidence only. Goldston's result does not imply the filtered variance formula needed here.

### 6.3 Open problem 1: the true remainder

Find a uniform estimate, in a stated regime of H, Q, k and X, for

$$\frac{1}{X}\int_X^{2X} \left| \omega_k \int_0^H e^{i\omega_k u} \left[S(T+u) - S_Q(T+u)\right] du \right|^2 dT,$$

and compare it with the diagonal term V of 4.2 without assuming that the two coincide. Three cases should be treated separately: modes far from omitted frequencies, modes near a resonance $\omega_k \approx \log q$, and sharp windows, whose Fourier transform decays only like the inverse of the frequency distance. A smooth window is a natural intermediate step. The formulas of Fujii and Goldston for the second moment of $S(t+h) - S(t)$ under RH are the natural starting point: the present problem replaces the difference by a windowed Fourier coefficient and subtracts $S_Q$.

### 6.4 Open problem 2: propagation to the coefficients

By 3.6 the error of the coefficients is the residual sequence $b_j$ convolved with the coefficients of $(1 - w)^{-\Delta S}$, and the amplification is bounded explicitly. Two steps remain. The first is a bound for $B_K = \Sigma|b_j|^2$ in terms of the residual power sums $r_k$, uniform in the window. The second is an estimate of the variation term $\Sigma\, \Delta b_j\, h_{k-j}$ of 3.6 for k below $\kappa_Q$, which would turn $\Phi_K(\Delta S)$ into an estimate for positive $\Delta S$. With Open problem 1 this would give a mean-square estimate for the coefficients, conditional on $\Delta S$.

### 6.5 Further checks

- Certified zeros and other heights. Everything here is at $T \approx 10^6$.
- Windows running from a zero to a zero, where N is fixed but the endpoints are singular for $\log \zeta$.
- The pointwise formula of 4.1 in mean square, which needs smooth windows.

### 6.6 Related observations, not developed here

- **Dirichlet L-functions.** For the characters mod 4 and mod 3, using zeros of the main sum of the approximate functional equation, peaks appear with intensity $|\chi(q)|^2 / q$ and are absent where $\chi(q) = 0$. Those zeros are not certified zeros of L.
- **CMV chain.** For the counting measure on a window, the first Verblunsky coefficients show prime powers with $|\alpha|^2 \approx \Lambda(q)^2 / (q\, \tau^2)$, and the last ones show all integers. With random Dirichlet weights the same peaks survive. For Verblunsky coefficients in matrix models of circular ensembles see Killip and Nenciu.

- **Failure of the proposed Bernstein–Szegő construction.** The lower half of the window polynomial had on average 66 of its 173 zeros inside the unit disk, in all 288 windows tested. So the measure built from that polynomial, as in the commensurate model of 2.3, does not exist for zeta windows. Other canonical measures are not excluded.

### 6.7 Scope

Nothing here bears on the truth of the Riemann hypothesis. The contribution, if it is new, is a quantitative description of how truncated prime-power sums, the zero count and the known smooth term determine the secular coefficients of a window of zeros. The only analytic error left is the part of S(t) not captured by the model.

## Appendix. Numerical methods and limits

No computation in this note is certified. The zeros come from a Riemann–Siegel code written for this note by an AI model and all statistics are from one height.

### A.1 Zeros

- Z(t) by the Riemann–Siegel formula with correction terms C0–C2, double precision with phases reduced in long double.
- Sign changes on a grid of 24 points per mean spacing, a parabola test for close pairs, then bisection.
- Validation at $T = 10^4$: 2 534 zeros, largest difference from an independent list $2 \times 10^{-11}$.
- At $T = 10^6$: 100 215 zeros found, against 100 215.25 expected from the smooth count. Independent check with mpmath, which is not a certificate: $N(T_0) = 1\,747\,146$ at $T_0 = 1\,000\,000.5741$; the counts of all 289 and 72 windows equal the differences of nzeros at their endpoints; 300 zeros chosen at random agree with zetazero to $2.3 \times 10^{-10}$.

### A.2 Other ingredients

- Windows: $T_j = \gamma_1 - 0.01 + jH$ as in 5.1; zeros with $T_j < \gamma \le T_j + H$.
- $\vartheta(t)$ from the complex log-gamma function; $\Delta S$ from $\vartheta$ and the zero count.
- $A_k$ from two integrations by parts, with $\vartheta'$ from the digamma function; it agrees with Simpson quadrature to $10^{-9}$.
- Power sums and Newton's recurrence in double precision; checked against 60 digits on the worst window.
- Model sums over prime powers $q \le Q$ with weight $\Lambda(q) / (\sqrt{q} \log q)$.
- Diagonal term V: primes and prime powers up to $10^7$, plus the tail $\omega_k^2 / (2\pi^2 \log^2 R)$.
- Bootstrap: moving blocks of 6 consecutive windows (H = 181) or 4 (H = 724), percentile intervals; 4 000 resamples, 1 500 for correlations. Adjacent windows share an endpoint, so $\Delta S$ has lag-1 autocorrelation −0.5. Intervals from independent resampling are close to these.
- Check of 3.4: Fourier coefficients by FFT on $2^{21}$ midpoint nodes. Check of 3.5: Newton's recurrence up to k = 6 000 in six windows. Check of 3.6: $b_j$ from the convolution of $c_k - \tilde{c}_k$ with the coefficients of $(1 - w)^{\Delta S}$; convolving back reproduces the error to $4 \times 10^{-14}$.

- Check of 4.2: the integral over [X, 2X] in closed form, from all pairwise frequency differences.
- Figure 1 shows window $j = 3$ ($T = 1\,000\,543.5741$, $N = 345$, $\Delta S = 0.058$), chosen for a residual energy close to the median.

### A.3 Limits

- One height, $T \approx 10^6$; two window lengths; 289 and 72 windows.
- Fixed windows only. The exploratory chain experiments used windows from zero to zero.
- No interval arithmetic. Bootstrap intervals measure variability between windows only; the block bootstrap accounts for short-range dependence.
- The error term of 4.1 is quoted from the published lemma and was not re-derived.
- The inequality used in 4.2 is quoted from a later paper, not from the original article.
- The expansion over integers in 3.4 has no proved error term.
- For Dirichlet L-functions the main-sum zeros miss about 0.3% of true zeros in the zeta calibration, and windows of anomalous length were discarded.
- The literature search of 1.1 used web search, abstracts and excerpts. MathSciNet and zbMATH were not consulted, and the paper of Coram and Diaconis is known here only through Diaconis's survey.

## AI contributions and research methodology

This research was directed by the author, an experimental physicist who works on mathematical problems outside his main professional activity. Two large language models, GPT-6 (OpenAI) and Claude Fable 5.2 (Anthropic), made substantial contributions to it. They are not authors. The author takes responsibility for the mathematical claims, the numerical results, the citations and the text.

- **Configuration.** The models were used adversarially. The author passed the statements, proofs, numerical results and objections of each model to the other, chose which checks to run, and decided what was kept, revised or classed as open.
- **GPT-6.** Critical review of every draft, and checks of quoted sources. First formulations of several statements: the linear response of 2.4, the symmetric model of 3.1–3.2, the full predictor of 3.3, the effective exponent of 3.5, the convolution lemma of 3.6 and its form by parts. Counterexamples, among them the one in 3.6.
- **Claude Fable 5.2.** All computations: the zero-finding code, the experiments, the tests and the figures were programmed and run by the model in a sandboxed environment, at the author's request. The first observation of prime-power structure in the coefficients, the commensurate model of 2.3, the formula of 4.1, the identities of 3.4, the final form of the proofs of 3.5 and 4.2, the constant-residual law of 3.6, the literature search of 1.1, and the drafting of the text.
- **Corrections made in the process.** A claimed decay of the remainder in 3.4 was withdrawn after a wider test. The conjecture that $\delta_Q$ governs the measured errors was tested and not confirmed (5.7). A statement that a negative jump cannot amplify the residual was refuted by a counterexample (3.6). A reported near-saturation of Young's bound was traced to windows with ΔS near 0 (5.7).

Neither model is an independent certifier of correctness, and errors may remain. The computations have not been re-run outside the model's sandbox; the scripts and the list of zeros are available from the author. These acknowledgments do not imply that OpenAI or Anthropic endorse the paper or its conclusions.

## References

The bibliographic data of every entry were checked on 19 November 2026 against arXiv records, publisher or repository pages, or the reference lists of later papers. The two statements used in proofs, Lemma 1 of Ford–Zaharescu and the Montgomery–Vaughan inequality as quoted by Bourgade, were checked against those papers. The papers of Goldston, Fujii, Gonek, Montgomery–Vaughan and Coram–Diaconis were not consulted directly: they are cited for statements quoted in the later sources named in each entry. Landau's paper was not consulted either: it is cited for the classical formula, and its bibliographic data were checked on EuDML. The other items of Section 1.1 were read through abstracts and excerpts, and Tao's two posts in part.